\documentclass[12pt]{article}
\usepackage{amssymb,amsmath,amsthm, amsfonts}

\def\sqr#1#2{{\vcenter{\vbox{\hrule height.#2pt
              \hbox{\vrule width.#2pt height#1pt \kern#1pt \vrule width.#2pt}
          \hrule height.#2pt}}}}
\def\signed #1{{\unskip\nobreak\hfil\penalty50
          \hskip2em\hbox{}\nobreak\hfil#1
          \parfillskip=0pt \finalhyphendemerits=0 \par}}
\def\endpf{\signed {$\sqr69$}}
\def\dbR{\hbox{\rm l\negthinspace R}}

\def\d{\delta}

\def\cF{{\cal F}}

\def\sqr#1#2{{\vcenter{\vbox{\hrule height.#2pt
              \hbox{\vrule width.#2pt height#1pt \kern#1pt \vrule width.#2pt}
              \hrule height.#2pt}}}}
\def\signed #1{{\unskip\nobreak\hfil\penalty50
              \hskip2em\hbox{}\nobreak\hfil#1
              \parfillskip=0pt \finalhyphendemerits=0 \par}}
\def\endpf{\signed {$\sqr69$}}

\def\dbR{{\mathop{\rm l\negthinspace R}}}

\def\3n{\negthinspace \negthinspace \negthinspace }
\def\2n{\negthinspace \negthinspace }
\def\1n{\negthinspace }

\def\dbR{{\mathop{\rm l\negthinspace R}}}
\def\={\buildrel \triangle \over =}

\def\d{\delta}

\def\O{\Omega}

\def\cF{{\cal F}}

\def\ms{\medskip}
\def\bs{\bigskip}
\def\q{\quad}

\def\sup{\mathop{\rm sup}}

\def\inf{\hbox{\rm inf$\,$}}

\def\ae{\hbox{\rm a.e.{ }}}
\def\as{\hbox{\rm a.s.{ }}}

\def\|{\Big |}
\def\({\Big (}
\def\){\Big )}
\def\[{\Big[}
\def\]{\Big]}
\def\be{\begin{equation}}
\def\bel{\begin{equation}\label}
\def\ee{\end{equation}}
\def\bt{\begin{theorem}}
\def\bcd{\begin{condition}}
\def\ecd{\end{condition}}
\def\et{\end{theorem}}
\def\bc{\begin{corollary}}
\def\ec{\end{corollary}}
\def\bde{\begin{definition}}
\def\ede{\end{definition}}
\def\bl{\begin{lemma}}
\def\el{\end{lemma}}
\def\bp{\begin{proposition}}
\def\ep{\end{proposition}}
\def\br{\begin{remark}}
\def\er{\end{remark}}
\def\ba{\begin{array}}
\def\ea{\end{array}}
\def\ed{\end{document}}

\def\square#1{\vbox{\hrule\hbox{\vrule height#1%
     \kern#1\vrule}\hrule}}
\def\rectangle#1#2{\vbox{\hrule\hbox{\vrule height#1%
     \kern#2\vrule}\hrule}}

\font\tenbb=msbm10 \font\sevenbb=msbm7 \font\fivebb=msbm5

\newfam\bbfam
\scriptscriptfont\bbfam=\fivebb \textfont\bbfam=\tenbb
\scriptfont\bbfam=\sevenbb

\newtheorem{lemma}{Lemma}[section]
\newtheorem{remark}{Remark}[section]

\newtheorem{theorem}{Theorem}[section]
\newtheorem{corollary}{Corollary}[section]

\newtheorem{definition}{Definition}[section]
\newtheorem{proposition}{Proposition}[section]
\newtheorem{condition}{Condition}[section]

\makeatletter
   
   \@addtoreset{equation}{section}
\makeatother

\begin{document}

\title{ Local Strict Comparison Theorem and Converse Comparison Theorems for Reflected Backward Stochastic Differential Equations \bf\footnote{Partially
supported by the Natural Science Foundation of China under grants
10426022, 10325101 (distinguished youth foundation, entitled
with``Control Theory of Stochastic Systems"), the Chang Jiang
Scholars Programme, and the Science Foundation of Chinese Ministry
of Education under grant 20030246004. }}

 \author{Juan Li\thanks{Institute of Mathematics, School of Mathematical Sciences, Fudan University, Shanghai
200433, China;  Key Laboratory of Mathematics  for Nonlinear
Sciences (Fudan University), Ministry of Education; \& Department
of Mathematics, Shandong University at Weihai, Weihai 264200,
China. {\small\it E-mail:} {\small\tt juanli@sdu.edu.cn}.}  \quad
and \quad Shanjian Tang\thanks{Department of Finance and Control
Sciences, School of Mathematical Sciences, Fudan University,
Shanghai 200433, China; \& Key Laboratory of Mathematics  for
Nonlinear Sciences (Fudan University), Ministry of Education.
{\small\it E-mail:} {\small\tt sjtang@fudan.edu.cn}.\ms} }

\date{}

 \maketitle

\begin{abstract} A local strict comparison theorem and some converse
comparison theorems are proved for reflected backward stochastic
differential equations  under suitable conditions.
\end{abstract}



\section{Introduction}

\hskip1cm The comparison theorem  turns out to be one
   classical result  for backward stochastic
differential equations (BSDEs). It allows us to compare the
  solutions of two real-valued BSDEs by comparing the
  terminal conditions and the generators. In a converse way, Peng
  is concerned in 1997 with the following converse comparison property for BSDEs: if the
  solutions of two real-valued BSDEs are equal at the initial
  time for any identical terminal condition, their generators are
  identical. For works on this problem, the reader is referred to among others: Chen~\cite{C1}, Briand et al.~\cite{BCHMP},
  Coquet et al.~\cite{CHMP1},
 and  Jiang~\cite{J1}. In their arguments, the strict comparison theorem of BSDEs plays
a crucial role.

On the other hand, the solution $Y$ of a reflected BSDE (RBSDE)
characterizes the value process of an optimal stopping time
problem, and  the price process $\{Y_{t}\}_{0\le t\le T}$ of an
American option is a solution of an RBSDE (See El Karoui et al.
\cite{EPaQ}):
\begin{equation} Y_{t} =(X_{T}-k)^{+}-\int_{t}^{T}[rY_{s}+\theta
Z_{s}]ds+K_T-K_t-\int_{t}^{T} Z_{s}dB_{s},
\end{equation}
with
$$
Y_{t}\ge S_{t}:=(X_{t}-k)^{+},\q \forall t\in [0, T];\q
\int_0^T(Y_t-S_t)\, dK_t=0.
$$
Here $\theta := \sigma^{-1}(\mu -r)$ is the premium of the market
risk, and $\{X_t\}_{0\le t\le T}$ is the stock price process
satisfying the following SDE:
$$X_{t} =X_{0}+\int_{0}^{t}\mu X_{s}ds+\int_{0}^{t}\sigma X_{s}dB_{s}, \q t\in [0,T].$$
 Define the stopping time $\tau :=\inf \{t:
Y_t=S_t\},$ which is the time when the investor would take the
action to sell or buy the stock. The theory of RBSDEs existing in
the literature only reveals how the price  $Y_{t}$ depends on the
generator $g(y,z):=ry+\theta z, y\in \dbR, z\in \dbR$ and the
strike price $k$ (more generally speaking, the obstacle and the
terminal value) as well. It is natural to ask how the premium
$\theta$ (more generally speaking, the generator $g$) can be
obtained from a family of American options parameterized by the
strike price $k$? Then, the relation among the solution, the
generator and the obstacle becomes interesting.

In this paper, we are concerned with comparison theorems and
converse comparison theorems for RBSDEs under suitable conditions,
which reveal some monotonicity between the solution, and the
generator and the obstacle of a RBSDE.

The rest of the paper is organized as follows. In Section 2, we
 provide some preliminary results on BSDEs and RBSDEs.
 In Section 3, we first illustrate that quite different
from BSDEs, RBSDEs do not have the global comparison property.
Then we prove a local strict comparison theorem for RBSDEs. In
Sections 4 and 5, we discuss converse comparison properties for
RBSDEs when the obstacle is not previously given and when the
obstacle is previously given, respectively. Some interesting
comparison theorems are obtained in both cases.

  \section{Preliminaries}

\hskip1cm In this section, we give some basic results on BSDEs and
RBSDEs. They will be used in the subsequent sections.

Let $(\Omega,{\cal{F}}, P)$ be a  probability space and
$\{B_t\}_{t\geq0}$ be a d-dimensional standard Brownian motion
    on this space such that $B_0=0$. Denote by
    $\{{\cal{F}}_t\}_{t\geq0}$ the filtration generated by
    Brownian motion $\{B_t\}_{t\geq0}$: $ {\cal{F}}_t:=\sigma\{B_{s}, s\in[0,t]\} \vee{\cal{N}}, t \in [0, T],
    $ where $ {\cal{N}}$ is the set of all P-null subsets.
    Let $T > 0$ be a given real number.  For any positive integer n
    and $z\in {\mathbf{R}}^{n}$, $|z|$ denotes the Euclidean norm.

    Define the following two spaces of processes:
   $${\cal{H}}^{2}(0,T;{\mathbf{R}}^{n}):=\{\{\psi_t\}_{0\leq t\leq T}\ \mbox{ \rm is an $\mathbf{R}^n$-valued predictable
   process s.t. } E\int^T_0|\psi_t|^2\, dt<+\infty \} $$ and
     $$ {\cal{S}}^2(0, T; {\mathbf{R}}):=\{\{\psi_t\}_{0\leq t\leq T}\ \mbox{ \rm is a
    predictable process  s.t. } E[\sup\limits_{0\leq t\leq T}| \psi_{t} |^2]< +\infty
\}.$$

Consider function $g: \Omega\times[0,T]\times {\mathbf{R}} \times
{\mathbf{R}}^{d} \rightarrow {\mathbf{R}} $ such that  $\{g(t, y,
z)\}_{t\in [0, T]}$ is progressively measurable for each $(y,z)$
in ${\mathbf{R}} \times {\mathbf{R}}^{d}$. We make  the following
assumptions on $g $ throughout the paper.

\bs (A1) There exists a constant $K>0$  such that \as for $t\in
[0,T]$,
$$|g(t, y_{1}, z_{1}) - g(t, y_{2}, z_{2})|\leq K(|y_{1}-y_{2}| +
|z_{1}-z_{2}|), \forall t\in [0, T], \forall y_{1}, y_{2}\in
{\mathbf{R}}, z_{1}, z_{2}\in {\mathbf{R}}^d.$$

\ms (A2) The process $g(\cdot,0, 0) \in
{\cal{H}}^{2}(0,T;{\mathbf{R}}).$

\ms (A3) $ g(t,y, 0)= 0$ \as for any $  (t,y) \in [0, T]\times
{\mathbf{R}}.$

\ms (A4) The mapping $  t \mapsto  g(t, y, z)$ is continuous
$a.s.$ for any $ (y,z) \in {\mathbf{R}}\times {\mathbf{R}}^{d}.$

\bs

 \br  It is obvious that Assumption (A3) implies Assumption (A2).\er

 It is by now well known (see Pardoux and Peng~\cite{PaPe} for the proof) that under Assumptions (A1)
and (A2), for any random variable $\xi\in L^2(\O, \cF_T,P),$ the
BSDE \be y_t = \xi + \int_t^Tg(s,y_s,z_s)ds - \int^T_tz_s\,
dB_s,\q 0\le t\le T \label{BSDE}\ee
 has a unique adapted solution $(y^{T, g, \xi}, z^{T, g,
\xi})\in {\cal{S}}^2(0, T; {\mathbf{R}})\times
{\cal{H}}^{2}(0,T;{\mathbf{R}}^{d}).$

In the sequel, we always assume that $g $ satisfies (A1) and (A2).
We introduce the following operator $\varepsilon_{g,T}$: for any
$\xi\in L^2(\O, \cF_T,P)$, denote by $\varepsilon_{g,
   T}[\xi]$ and $\varepsilon_{g,T}[\xi|{\cal{F}}_{t}]$ the initial value $y^{T, g, \xi}_0 $ and the value $y^{T, g, \xi}_t $
   at time $t$ of the solution to BSDE~(\ref{BSDE}), respectively. For a stopping
   time $\tau$, the operator $\varepsilon_{g,\tau}$ can be defined
   in an identical way.

We give some basic results of BSDEs, including Lemmas 2.1, 2.2 and
2.3, which  can be found in Briand et al.~\cite{BCHMP} or
Peng~\cite{Pe1}, El Karoui et al.~\cite{ElPeQu}, and Jiang
~\cite{J1}, respectively.

\bl (Comparison Theorem) \label{L1} Assume that two fields $g_1$
and $g_2$ satisfy (A1) and (A2). Consider $ \xi_1,\ \xi_2 \in
L^{2}(\Omega, {\cal{F}}_{T}, P)$. We have

{\rm (i) }(Monotonicity) If \ $ \xi_1 \geq \xi_2$  and \ $ g_1
\geq g_2 \ a.s.$, then $\varepsilon_{g_1,T}[\xi_1] \geq
\varepsilon_{g_2,T}[\xi_2]$, and\ $
\varepsilon_{g_1,T}[\xi_1|{\cal{F}}_{t}]\geq
                 \varepsilon_{g_2,T}[\xi_2|{\cal{F}}_{t}] \ a.s.\ \hbox{\it for } t\in [0,T].$

{\rm (ii)}(Strict Monotonicity) If  $ \xi_1 \geq \xi_2, g_1 \geq
g_2$, and \  $P(\{\xi_1
> \xi_2\})> 0$, then $P(\{\varepsilon_{g,T}[\xi_1|{\cal{F}}_{t}]
> \varepsilon_{g,T}[\xi_2|{\cal{F}}_{t}]\})>0  \hbox{ \it for } t\in [0,T].$ In particular,  $
\varepsilon_{g,T}[\xi_1]
> \varepsilon_{g,T}[\xi_2].$
\el

 \bl Assume that the field $g$ satisfies
assumptions (A1) and (A2).  Consider the stopping time $\tau \leq
T$ and $ \xi\in L^{2}(\Omega, {\cal{F}}_{\tau}, P)$. Define
$$
  \overline{g}(t, y, z) := g(t, y, z)1_{[0, \tau]}(t),\ \  \forall (t,y,z)\in [0, T]\times {\mathbf{R}}\times {\mathbf{R}}^d.  $$
  Then $$\varepsilon_{g, \tau}[\xi] = \varepsilon_{\overline{g},
  T}[\xi] \hbox{ \it and }
   \varepsilon_{g, \tau}[\xi|{\cal{F}}_{t}] = \varepsilon_{\overline{g}, T}[\xi|{\cal{F}}_{t}]\ a.s. \hbox{ \it for } t \in [0, \tau].$$
\el

\bl Assume that two functions $g_{1}$ and $g_{2}$ satisfy
assumptions (A1), (A2) and (A4). Then the following two assertions
are equivalent:

 {\rm(i)} $\ \varepsilon_{g_{1},
\tau}[\xi] = \varepsilon_{g_{2}, \tau}[\xi]$ for each stopping
time $\tau\leq T$ and any $\xi \in L^{2}(\Omega, {\cal{F}}_{\tau},
   P).$

 {\rm(ii)} $\ g_{1}(t,y,z) = g_{2}(t,y,z)\q \as$ for any $(t, y,
z)\in [0, T]\times {\mathbf{R}}\times {\mathbf{R}}^d$.
  \el

   We introduce the conditional $g$-expectation (see Peng~\cite{Pe1,Pe2}, Chen~\cite{C1}, Coquet et al.~\cite{CHMP1,CHMP2}).
   Suppose $g$ satisfies (A1), (A3)
   and (A4). We set, for any stopping time $\tau$ taking values in
   $[0,T]$, $$\varepsilon_{g}[\xi|{\cal{F}}_{\tau}]:= y^{T, g, \xi}_{\tau}(=\varepsilon_{g,T}[\xi|{\cal{F}}_{\tau}]).$$
   It can be shown that $\varepsilon_{g}[\xi|{\cal{F}}_{\tau}]$ is the unique  ${\cal{F}}_{\tau}$-measurable, square-integrable random variable $\eta$ such
   that
   $$ \varepsilon_{g}[1_{A}\eta] =
   \varepsilon_{g}[1_{A}\xi], \quad \forall A \in {\cal{F}}_{\tau}.$$
   Therefore it is called the $g$-expectation conditioned on
   $\cF_\tau$.
      Notice that $g$-expectation $\varepsilon_{g}$ is a particular
   example of the nonlinear expectation introduced in~\cite{C1,CHMP2,Pe1,Pe2}.
    Now we borrow from~\cite{CHMP1} the converse comparison theorem for $g$-expectation.

\bl Suppose that two functions $g_{1}$ and $g_{2}$ satisfy
assumptions (A1), (A3) and (A4). Then the following two assertions
are equivalent:

{\rm(i)} $ \varepsilon_{g_{1}}[\xi] \geq \varepsilon_{g_{2}}[\xi]$
for any $\xi \in L^{2}(\Omega, {\cal{F}}_{T},
   P).$

{\rm(ii)} $\ g_{1}(t,y,z) \geq g_{2}(t,y,z)\q \as$ for any $ (t,
y, z)\in [0, T]\times {\mathbf{R}}\times {\mathbf{R}}^d$.
  \el

 A reflected BSDE is associated with a terminal
   condition $\xi \in L^{2}(\Omega,{\cal{F}}_{T}, P)$, a
   generator $g$, and an ``obstacle" process $\{S_t\}_{0\leq t \leq
   T}$. We make the following assumption:

\bs
 (A5)\ $\{S_t\}_{0\leq t \leq T}$ is a continuous process such that $\{S_t\}_{0\leq t \leq T}\in {\cal{S}}^2(0, T;{\mathbf{R}})$.
\bs

The solution of a RBSDE is a triple $(Y, Z, K)$ of
${\cF}_t$-progressively measurable processes taking values in
$\mathbf{R}\times\mathbf{R}^d\times\mathbf{R}_+$ and
   satisfying

 {\rm (i)} $Z \in {\cal{H}}^{2}(0,T;{\mathbf{R}}^{d}),\  \mbox{Y} \in {\cal{S}}^2(0, T; {\mathbf{R}}),$ and
 $ K_{T} \in L^{2}(\Omega,{\cal{F}}_{T}, P);$
\be \mbox{\rm (ii)}  \ Y_t = \xi + \int_t^Tg(s,Y_s,Z_s)ds + K_{T}
- K_{t} - \int^T_tZ_sdB_s,\quad t\in
[0,T];\qquad\qquad\label{RBSDE}\ee

 {\rm (iii)} $Y_t \geq S_t$\q a.s. for any $ t\in [0,T];$

{\rm (iv)} $\{K_{t}\}$ is continuous and increasing, $K_{0}=0$ and
$ \int_0^T(Y_t - S_t)dK_{t}=0.$ \vskip0.5cm

The following two lemmas are borrowed from El Karoui et
al.~\cite{EKPPQ}.

 \bl Assume that $g$ satisfies (A1) and (A2), $ \xi \in L^{2}(\Omega, {\cal{F}}_{T}, P)$, $\{S_t\}_{0\leq t \leq T}$
satisfies (A5), and $S_T \leq \xi\ \ a.s..$ Then
RBSDE~(\ref{RBSDE}) has a unique solution $(Y, Z, K).$\el

\br For simplicity, a given triple $(\xi, g, S)$ is said to
satisfy the Standard Assumptions if the generator $g$ satisfies
(A1) and (A2), the terminal value $ \xi \in L^{2}(\Omega,
{\cal{F}}_{T}, P)$, the obstacle $S$ satisfies (A5), and $S_T \leq
\xi \ a.s..$ \er

\bl (comparison theorem) Suppose that two triples $(\xi_1, g_1,
S^1)$ and $(\xi_2, g_2, S^2)$ satisfy the Standard Assumptions (in
fact, it is sufficient for either $g_1$ or $ g_2$ to satisfy the
Lipschitz condition (A1)). Furthermore, we make the following
assumptions:
$$
  \begin{array}{ll}
{\rm(i)}&\xi_1 \leq \xi_2\ \ a.s.;\\
{\rm(ii)}&g_1(t,y,z) \leq g_2(t,y,z)\  \as \hbox{ \it for } (t,y,z)\in [0,T]\times {\mathbf{R}}\times {\mathbf{R}}^{d};\\
{\rm(iii)}& S_t^1 \leq S^2_t\ \ a.s. \hbox{ \it for } t\in [0,T]. \\
 \end{array}
  $$
  Let $(Y^1,Z^1, K^1)$ and $(Y^2, Z^2, K^2)$ be adapted solutions of RBSDEs~(\ref{RBSDE}) with data $(\xi_1, g_1,
  S^1)$ and $(\xi_2, g_2, S^2),$ respectively.  Then $Y^1_{t} \leq Y^2_{t}\ a.s. $ for $t\in [0,T].$\el

\bl \label{L7} Assume that $(\xi_{1}, g_{1}, S)$ and $(\xi_{2},
g_{2}, S)$ satisfy the Standard Assumptions. Furthermore, we make
the following assumptions:
$$
  \begin{array}{ll}
{\rm(i)}&\xi_{1} \leq \xi_{2}\ a.s.;\\
{\rm(ii)}&g_{1}(t,y,z) \leq g_{2}(t,y,z)\ \as \hbox{ \it for } (t,y,z)\in [0,T]\times {\mathbf{R}}\times {\mathbf{R}}^{d}.\\
 \end{array}
  $$
  Let $(Y^{1},Z^{1}, K^{1})$   and $(Y^{2},Z^{2}, K^{2})$ be  adapted solutions of RBSDEs~(\ref{RBSDE}) with data $(\xi_{1}, g_{1},
  S)$ and $(\xi_{2}, g_{2}, S)$, respectively.
  Then we have
  \be K^{1}_{t} \geq K^{2}_{t} \ \ a.s. \hbox{ \it for } t\in [0, T],  \mbox{ \it and }  K^{1}_{t}- K^{2}_{t} \ \mbox{ \it is increasing
    in time variable}\  t.\label{Property}\ee
    \el

See Hamad\`ene et al.~\cite[Proposition 41.3]{HLM} for the
detailed proof of Lemma~\ref{L7}.

\section{Local strict comparison theorem for RBSDEs }

 \hskip1cm In contrast to BSDEs, the strict comparison theorem is not true
  in general for RBSDEs. Here are two counterexamples.

 {\bf Example 3.1}\ Take $T=1, \ g =
\frac{1}{3},\  \xi_{1}= \frac{1}{3},\ \xi_{2}= \frac{1}{2}$,\ and
$ S_{t}=-2t+1$ for $t\in [0, T]$. Then the solution $(Y^{1}
,Z^{1},K^{1})$ of RBSDE~(\ref{RBSDE}) with data $(\xi_{1}, g, S)$
is given by
$$\begin{aligned}
  Y^{1}_t=
   \left\{
  \begin{array}{ll}
  1-2t,&\mbox{if } 0\leq t \leq \frac{1}{5}\\
   \frac{2}{3}-\frac{1}{3}t,&\mbox{if }\frac{1}{5}< t \leq 1;
  \end{array}
  \right.\\
  \end{aligned}
\quad Z^{1}_t = 0,\ \ 0\leq t\leq 1;\quad
\begin{gathered}
  K^{1}_t=
   \left\{
  \begin{array}{ll}
   \frac{5}{3}t,
      &\mbox{if } 0\leq t \leq \frac{1}{5}\\
   \frac{1}{3},&\mbox{if }\frac{1}{5}< t \leq
   1.
  \end{array}
  \right.\\
  \end{gathered} $$
The solution $(Y^{2}, Z^{2}, K^{2})$ of RBSDE~(\ref{RBSDE}) with
data $(\xi_{2}, g, S)$ is given by
  $$\begin{aligned}
  Y^{2}_t=
   \left\{
  \begin{array}{ll}
  1-2t, &\mbox{if } 0\leq t \leq \frac{1}{10}\\
   \frac{5}{6}-\frac{1}{3}t,&\mbox{if }\frac{1}{10}< t \leq 1;
  \end{array}
  \right.\\
   \end{aligned}
\quad Z^{2}_t = 0,\ \ 0\leq t\leq 1;\quad
\begin{gathered}
  K^{2}_t=
   \left\{
  \begin{array}{ll}
   \frac{5}{3}t, &\mbox{if } 0\leq t \leq \frac{1}{10}\\
   \frac{1}{6},&\mbox{if }\frac{1}{10}< t \leq 1.
  \end{array}
  \right.\\
   \end{gathered}$$

Obviously, $Y^{1}_t=Y^{2}_t$ for $ t\in [0, \frac{1}{10}]$ and $
Y^{1}_t < Y^{2}_t$ for $t\in (\frac{1}{10}, 1]$. Moreover,
$K^{1}_t\geq K^{2}_t$ for $t \in [0,1].$ The strict comparison
theorem of RBSDE~(\ref{RBSDE}) does not hold.

The following example shows that even if the generator is zero, it
happens that the strict comparison theorem of RBSDE~(\ref{RBSDE})
may not be true.

 {\bf Example 3.2}\ Take $T=1, \ g = 0,\
\xi_{1}= \frac{1}{3},\ \xi_{2}= \frac{1}{2}$, and $S_{t}=-2t+1$
for $t\in [0, T]$. Then the solution $(Y^{1}, Z^{1}, K^{1})$ of
RBSDE~(\ref{RBSDE}) with data $(\xi_{1}, g, S)$ is given by
$$\begin{aligned}
  Y^{1}_t=
   \left\{
  \begin{array}{ll}
  1-2t,&\mbox{if } 0\leq t \leq \frac{1}{3}\\
   \frac{1}{3},&\mbox{if }\frac{1}{3}< t \leq 1;
  \end{array}
  \right.\\
  \end{aligned}
\quad Z^{1}_t = 0,\ \ 0\leq t\leq 1;\quad
\begin{gathered}
  K^{1}_t=
   \left\{
  \begin{array}{ll}
   2t,&\mbox{if } 0\leq t \leq \frac{1}{3}\\
   \frac{2}{3},&\mbox{if }\frac{1}{3}< t \leq 1.
  \end{array}
  \right.\\
  \end{gathered} $$
The solution $(Y^{2}, Z^{2}, K^{2})$ of RBSDE~(\ref{RBSDE}) with
data $(\xi_{2}, g, S)$ is given by
  $$\begin{aligned}
  Y^{2}_t=
   \left\{
  \begin{array}{ll}
  1-2t, & \mbox{if } 0\leq t \leq \frac{1}{4}\\
   \frac{1}{2},& \mbox{if }\frac{1}{4}< t \leq 1;
  \end{array}
  \right.\\
   \end{aligned}
\quad Z^{2}_t = 0,\ \ 0\leq t\leq 1;\quad
\begin{gathered}
  K^{2}_t=
   \left\{
  \begin{array}{ll}
   2t,&\mbox{if } 0\leq t \leq \frac{1}{4}\\
   \frac{1}{2},&\mbox{if }\frac{1}{4}< t \leq
   1.
  \end{array}
  \right.\\
   \end{gathered}$$

Obviously, $Y^{1}_t=Y^{2}_t$ for $ t\in [0, \frac{1}{4}]$ and $
Y^{1}_t < Y^{2}_t$ for $t\in (\frac{1}{4}, 1].$ Moreover,
$K^{1}_t\geq K^{2}_t$ for $t \in [0,1].$ The strict comparison
theorem of RBSDE~(\ref{RBSDE}) does not hold.

However, we have the local strict comparison theorem.

\bt \label{T1} Suppose that two triples $(\xi_{1}, g, S)$ and
$(\xi_{2}, g, S)$ satisfy the Standard Assumptions. Moreover,
assume that
$$\xi_{1} \leq \xi_{2} \q \as \mbox{ \it and }  P(\{\xi_{1} < \xi_{2}\})>0. $$
  Let $(Y^{1},Z^{1}, K^{1})$ and $(Y^{2},Z^{2}, K^{2})$ be adapted solutions of RBSDE~(\ref{RBSDE}) with data $(\xi_{1}, g,
  S)$ and $(\xi_{2}, g, S)$, respectively.
   Then there exists a stopping time $\tau$ such that
    $\tau < T$ almost surely and $P(\{Y_t^1<Y_t^2,\  \forall t\in
[\tau, T]\})>0$.    \et

   {\bf Proof.} From Lemma 2.6, we have  $$Y^{1}_{t} \leq Y^{2}_{t} \  \ a.s. \hbox{ \it for } t\in [0, T].\eqno(3.1)
    $$

    Define a sequence of stopping times $\{\tau_k\}_{k=1}^\infty$ in the following way:
  $$\tau_1:=0$$
  and
    $$\tau_{k+1}= \mbox{inf}\{t\geq \tau_k+{1\over 2}(T-\tau_k): Y^{1}_{t} = Y^{2}_{t}\}\wedge T, \q k=1,2,\ldots.$$
    It is obvious that the sequence $\{\tau_k\}_{k=1}^\infty$ is both
    bounded by $T$ and non-decreasing. Therefore, it has an almost sure limit $\tau$,
    which is still a stopping time satisfying $\tau\le T$. Since
$$
\tau_{k+1}\geq \tau_k+{1\over 2}(T-\tau_k), \q  k=1,2,\ldots,
$$
we have by passing to the limit that $\tau\ge \tau+ {1\over
2}(T-\tau)$, that is, $\tau\ge T.$ Hence, $\tau=T$.

Furthermore, we assert that there is some positive integer $k_0$
such that $P(\{\tau_{k_0}=T\})>0$. Otherwise, we have
    $$
P(\{\tau_k<T\})=1, \q  k=1,2,\ldots.
    $$
This implies the following
$$
Y^1_{\tau_k}=Y^2_{\tau_k}, \q  k=1,2,\ldots.
$$
Passing to the limit, we have $Y^1_T=Y^2_T\ \as$. That is,
$\xi_1=\xi_2$, which contradicts the assumption that $P(\{\xi_{1}
< \xi_{2}\})>0$.

Take the smallest integer $\tilde k$ among those positive integers
$k_0$ such that $P(\{\tau_{k_0}=T\})>0$. Then, we have
$$\tilde
k\ge 1,\q  P(\{\tau_{\tilde k}=T\})>0, \q \tau_{\tilde k-1}<T\
\as.$$
 We assert that the stopping
time
$$\tilde \tau:=[\tau_{\tilde k-1}+{1\over
2}(T-\tau_{\tilde k-1})]<T$$ is a desired one of the theorem. In
fact, by definition of $\tau_{\tilde k}$, we have $Y_t^1<Y_t^2$
 on the interval $[\tilde\tau, T]$ whenever
$\tau_{\tilde k}=T$. Therefore, \be P(\{Y_t^1<Y_t^2, \  \forall
t\in [\tilde\tau, T]\})\ge P(\{\tau_{\tilde k}=T\})>0.  \ee The
proof is complete.\endpf

\bc In Theorem~\ref{T1}, if $g$ is either bounded from below by a
nonnegative constant $C_1$ or satisfies (A3), $\ \xi_{i}$ ($i = 1,
2$) are bounded from below and the obstacle process $S$ are
bounded from above by a constant $C_2$, then we have
$$Y^{1}_t \leq Y^{2}_t \  a.s., \mbox{\it and } P(\{Y^{1}_t < Y^{2}_t\})>0 \hbox{ \it for }  t\in [0, T].
$$
In particular, we have $Y^{1}_0 < Y^{2}_0$.\ec

{\bf Proof.} Consider the following BSDEs:
$$Y'^{i}_t = \xi_{i} + \int_t^{T}g(s, Y'^{i}_s, Z'^{i}_s)ds - \int^T_tZ'^{i}_sdB_s,\ \ \  0\leq t\leq
T,\ \ i=1, 2.$$ Obviously, it follows from Lemma~\ref{L1} that
$Y'^{i}_t \geq C_2 \geq S_{t} \ a.s., i=1,2, $ and
$P(\{Y_t'^1<Y_t'^2\})>0$ for $t\in [0, T].$ From Lemma 2.5, we
have
 $$ Y^{i}_t = Y'^{i}_t,\  Z^{i}_t = Z'^{i}_t,\hbox { \it and } K^{i}_t= 0 \ \   a.s. \hbox{ \it for } t\in [0, T].$$
 Therefore  $ P(\{Y^{1}_t < Y^{2}_t\})>0$ for $t\in [0,
T].$  Then  the proof is complete. \endpf

\section{ A converse problem for RBSDEs }

\hskip1cm In this section, we consider the general converse
comparison theorem for RBSDE~(\ref{RBSDE}).  Coquet et
al.~\cite{CHMP1} prove the converse comparison for BSDEs: if $g_1,
g_2$ satisfy (A1), (A3) and (A4), and $\varepsilon_{g_{1}, T}[\xi]
\geq \varepsilon_{g_{2}, T}[\xi]$ for any $\xi \in L^{2}(\Omega,
{\cal{F}}_{T}, P)$,  then $g_{1}(t, y, z) \geq g_2(t, y, z)\ a.s.$
for $(t,y,z)\in [0,T]\times \mathbf{R}\times \mathbf{R}^d$.

Consider the converse comparison  for RBSDE~(\ref{RBSDE}). It is
interesting since   the strict comparison theorem is not true for
RBSDE~(\ref{RBSDE}) and several arguments developed for BSDEs have
to be modified for RBSDEs.

 In the sequel, we always assume that the data $(\xi, g, S) $ satisfies the Standard Assumptions for RBSDEs.
 We introduce the following operator $\varepsilon_{g}^{r}$: denote by $\varepsilon^{r}_{g,
   T}[\xi]$ and
   $\varepsilon^{r}_{g,T}[\xi|{\cal{F}}_{t}]$ the initial value $Y^{T, g, \xi, S}_0 $
   and the value $Y^{T, g, \xi, S}_t $ at time $t$ of the solution of RBSDE~(\ref{RBSDE}) with data $(\xi, g, S)$, respectively.

\bt \label{T2} Suppose that two functions $g_{1}$ and $g_{2}$
satisfy assumptions (A1), (A3) and (A4). Then the following two
conditions are equivalent:

{\rm (i)} $\varepsilon^{r}_{g_{1}, T}[\xi] \geq \varepsilon^{r}_{g_{2},
   T}[\xi]$ for any $ \xi \in L^{2}(\Omega, {\cal{F}}_{T},
   P)$ and any obstacle process $(S_t)_{0\leq t \leq T}$ satisfying
   (A5) and $\xi \geq S_{T}\ \ a.s.$.

{\rm (ii)} $g_{1}(t,y,z) \geq g_{2}(t,y,z) \ \as$ for  $(t, y,
z)\in [0,T]\times \mathbf{R}\times \mathbf{R}^d$.

\et

{\bf Proof.} Thanks to Lemma 2.6, it is obvious that (ii) implies
(i). It is sufficient to prove that (i) implies (ii).

 For $\xi \in L^{2}(\Omega, {\cal{F}}_{T}, P)$,
 consider the following BSDE:
$$
  \left\{
  \begin{array}{ll}
   -dy^{i}(t) = g_{i}(t, y^{i}(t), z^{i}(t))dt - z^{i}(t)dB_{t},\q t\in [0,T],\\
   y^{i}(T)= \xi,\\
  \end{array}
  \right.\\
   $$
  for $i=1, 2.$
In view of Lemma 2.4, it suffices to show that
$\varepsilon_{g_{1}, T}[\xi] \geq  \varepsilon_{g_{2}, T}[\xi].$

Consider the following BSDE:
   $$Y_t =  \xi + \int^{T}_{t}[- K|Y_s| - K|Z_{s}|]ds - \int^{T}_{t}Z_{s}dB_s,\q 0 \leq t \leq T,$$
where $K=\mbox{max}(K_1, K_2)$ with $K_1$ and $K_2$ being the
Lipschitz constants of $g_1$ and $g_2$, respectively.
   Then $\{Y_t\}_{0\leq t \leq T}$  satisfies (A5) (see El Karoui et al.~\cite{ElPeQu}  for the detailed proof). Set $S:=Y$.
   Since $\xi\geq S_{T}$, from assumptions (A1), (A3) and Lemma 2.1, it follows that
   $$y^{i}(t)\geq S_{t} \q \as \hbox{ \it for }\  t\in [0,T].$$
   From Lemma 2.5, we see that $(y^{i}, z^{i}, 0)$ is the solution of RBSDE~(\ref{RBSDE}) with data
$(\xi, g_i, S)$ for $i=1,2$.  In particular,
$\varepsilon^{r}_{g_{i}, T}[\xi] = y^{i}(0)$ for $i=1, 2.$
Therefore, we get from the assumption that
$$\varepsilon_{g_{1}, T}[\xi] = \varepsilon^{r}_{g_{1}, T}[\xi] \geq \varepsilon^{r}_{g_{2}, T}[\xi] = \varepsilon_{g_{2}, T}[\xi] .$$
The proof is complete. \endpf

As an immediate consequence of Theorem~\ref{T2}, we have

\bc Assume that functions $g_{1}$ and $g_{2}$ satisfy assumptions
(A1), (A3) and (A4). Then the following two conditions are
equivalent:

{\rm (i)} $\varepsilon^{r}_{g_{1}, T}[\xi] =
\varepsilon^{r}_{g_{2}, T}[\xi]$ for any $\xi \in L^{2}(\Omega,
{\cal{F}}_{T}, P)$ and any obstacle process $\{S_t\}_{0\leq t \leq
T}$ satisfying (A5) and $\xi \geq S_{T} \ a.s.$.

{\rm (ii)} $ g_{1}(t,y,z) = g_{2}(t,y,z)\ a.s.$ for any $(t, y,
z)\in [0,T]\times\mathbf{R}\times \mathbf{R}^d$. \ec

\br  The assertion of Theorem~\ref{T2} is not true if assumption
(A3) fails to be satisfied. To show this fact, consider the
following
 example:
  $$g_{1}(t)= t 1_{[0,\frac{T}{2})}(t) + \frac{T}{2} 1_{[\frac{T}{2},
  T]}(t) \hbox{ \it and }  g_{2}(t)= \frac{T}{2} 1_{[0,\frac{T}{2})}(t) + (T-t) 1_{[\frac{T}{2},
  T]}(t).
  $$
  Then, for any $\xi \in L^{2}(\Omega, {\cal{F}}_{T},
   P)$,\ any $\{S_t\}_{0\leq t \leq T}$ satisfying (A5), and $\xi \geq S_{T}
  \ a.s.$, it follows from
   El Karoui et al.~\cite[Proposition 2.3]{EPaQ} that
   $$\varepsilon^{r}_{g_{1}, T}[\xi] \leq \varepsilon^{r}_{g_{2}, T}[\xi].$$
    However,   $g_{1}\not\leq g_{2}.$\er

\br If the obstacle process $\{S_t\}_{0\leq t \leq T}$ is
previously given and $\varepsilon^{r}_{g_{1}, T}[\xi] \geq
\varepsilon^{r}_{g_{2}, T}[\xi]$ only  for those $\xi \in
L^{2}(\Omega, {\cal{F}}_{T}, P)$ such that $\xi \geq S_{T}\ a.s.$,
then Theorem~\ref{T2} is not true in general. It suffices to
consider the following example:
$$g_1(t, y, z)= \mu_1(y-c_1)^{-}\wedge|z| \hbox{ \it and }  g_2(t, y, z)= \mu_2(y-c_2)^{-}\wedge|z|,$$
with $c_1< c_2$ and $\mu_1 > \mu_2$. It is obvious that $g_1$ and
$g_2$ satisfy assumptions (A1), (A3), and (A4). Furthermore, take
$S_t = c_2,\ t\in [0,T]$. Then for any $\xi \in L^{2}(\Omega,
{\cal{F}}_{T}, P)$ satisfying $\xi \geq S_{T}\ a.s.$,  we have
$$\varepsilon^{r}_{g_{1}, T}[\xi] = \varepsilon^{r}_{g_{2},
T}[\xi] \hbox{ \it and}\
 \varepsilon^{r}_{g_{1}, T}[\xi|{\cal{F}}_{t}] = \varepsilon^{r}_{g_{2}, T}[\xi|{\cal{F}}_{t}]\ \ a.s. \hbox{ \it for } t\in (0,
    T].$$
However,  we  have neither $g_1\geq g_2$ nor  $ g_1 \leq g_2$.\er

 When assumption (A2) instead of (A3) is made on $g$, we have the following converse comparison result.

\bt \label{T3} Assume that functions $g_{1}$ and $g_{2}$ satisfy
assumptions (A1), (A2) and (A4). Then the following two conditions
are equivalent:

{\rm (i)} $\varepsilon^{r}_{g_{1}, \tau}[\xi] =
\varepsilon^{r}_{g_{2}, \tau}[\xi]$ for any stopping time $\tau
\leq T$, any $\xi \in L^{2}(\Omega, {\cal{F}}_{\tau}, P)$, and any
obstacle $\{S_t\}_{0\leq t \leq T}$ satisfying (A5) and $\xi\geq
S_{\tau}\  a.s..$

{\rm (ii)}  $ g_{1}(t,y,z) = g_{2}(t,y,z)$ $a.s.$ for any $(t, y,
z)\in [0,T]\times\mathbf{R}\times \mathbf{R}^d$. \et

 {\bf Proof.} Thanks to Lemma 2.6, it is obvious that (ii) implies (i). It is
sufficient to  prove that (i) implies (ii).

 For $\xi \in L^{2}(\Omega, {\cal{F}}_{\tau}, P)$,
 consider the following BSDE defined on the interval $[0, \tau]:$
$$
  \left\{
  \begin{array}{ll}
   -dy^{i}(t) = g_{i}(t, y^{i}(t), z^{i}(t))dt - z^{i}(t)dB_{t},\\
   y^{i}(\tau)= \xi,\\
  \end{array}
  \right.\\
   $$ for $i=1, 2$.
In view of Lemma 2.3, it suffices to show that
$\varepsilon_{g_{1}, \tau}[\xi] =  \varepsilon_{g_{2},
\tau}[\xi].$

  Consider the obstacle process  $\{S_{t}\}_{0\leq t \leq T}$ which is defined to be $\xi$ on
   $[\tau, T]$, and on $[0, \tau]$ is taken to be one component of the solution of the following BSDE:
   $$S_t = \xi + \int^{\tau}_t[g_1(s, 0, 0)\wedge g_2(s, 0, 0) - K|S_s| - K|Z_{s}|]ds - \int^{\tau}_tZ_{s}dB_s,$$
   where $K=\mbox{max}(K_1, K_2)$, $K_1$ and $K_2$ are the Lipschitz constants of $g_1$ and $g_2$, respectively. The above equation admits a unique
solution $\{S_{t}\}_{0\leq t \leq T}$, which satisfies (A5). Since
$\xi\geq S_{\tau}$, it follows from assumptions (A1), (A2), Lemmas
2.1 and 2.2 that
   $$y^{i}(t)\geq S_{t}\ a.s.\  \mbox{ on the interval}\ [0, \tau].$$
   From Lemma 2.5, we get that $(y^{i}, z^{i}, 0)$ is the solution of RBSDE~(\ref{RBSDE}) with data
$(\xi, g_i, S)$\ on the interval $[0, \tau]$  for $i=1,2$.  In
particular, $\varepsilon^{r}_{g_{i}, T}[\xi] = y^{i}(0)$ for $i=1,
2.$ Therefore, it follows from the assumption that
$$\varepsilon_{g_{1}, \tau}[\xi] = \varepsilon^{r}_{g_{1}, \tau}[\xi] = \varepsilon^{r}_{g_{2}, \tau}[\xi] = \varepsilon_{g_{2}, \tau}[\xi] .$$
The proof is complete.\endpf

\br From the example in Remark 4.1, we can get that if $g$ only
satisfies (A1), (A2), and (A4), assertion {\rm(ii)} of
Theorem~\ref{T3} fails to hold under the condition {\rm(i)} of
   Corollary 4.1. \er

\section{Alternative  converse problem for RBSDEs with the obstacle process
$\{S_t\}_{0\leq t \leq T}$ being given}

 \hskip1cm Remark 4.2 shows that if the obstacle process $\{S_t\}_{0\leq t \leq T}$ is previously given, it is impossible in general to
   compare  the generator $g$ on the whole space
 $\O\times [0,T]\times \mathbf{R}\times \mathbf{R^d}$. In this section we shall show that  we
can still have the local converse comparison theorem for RBSDEs on
an upper semi-space $\Omega\times [0,T]\times [C,+\infty)\times
\mathbf{R}^d$, specified by the uniform upper bound $C$ of the
obstacle, which is actually the whole space if the generator does
not depend on the first unknown variable $y$ (see Theorem~\ref{T5}
below).

Assume that the data $(\xi, g, S)$ satisfies the Standard
Assumption for RBSDEs. But to emphasize the dependence on the
obstacle process $\{S_t\}_{0\leq t \leq T}$, denote by
$\varepsilon^{r, S}_{g,
   T}[\xi]$ and $\varepsilon^{r, S}_{g,T}[\xi|{\cal{F}}_{t}]$ the initial value $Y^{T, g, \xi, S}_0 $ and  the  value $Y^{T, g, \xi, S}_t $
at time $t$ of the solution of RBSDE (2.2) with data  $(\xi, g,
S)$, respectively.

\bp Suppose that $g$ satisfies (A1) and (A2), and the obstacle
process $\{S_t\}_{0\leq t \leq T}$ satisfies (A5). For the
stopping time $\tau \leq T$, and the terminal value $\xi\in
L^{2}(\Omega, {\cal{F}}_{\tau}, P)$ such that $\xi \geq S_{\tau}\
a.s.$, then we have
$$\varepsilon^{r, S}_{g,\tau}[\xi] = \varepsilon^{r,
\overline{S}}_{\overline{g},T}[\xi]$$ where
 $$ \overline{g}(t, y, z) :=  g(t, y, z)1_{[0, \tau]}(t) \hbox{ \it and }  \overline{S}_{t} :=
 S_{t\wedge\tau} \ \as \hbox{ \it for }
  (t, y, z)\in [0,\ T]\times {\mathbf{R}}\times {\mathbf{R}}^{d}.$$
\ep

{\bf Proof.} Consider the solution
 $(Y^{\tau, g, \xi, S}, Z^{\tau, g, \xi, S}, K^{\tau, g, \xi, S})$ of RBSDE (2.2)
  with data $(\xi, g, S) $ on the interval $[0, \tau]$, and
 $(\overline{Y}^{T, \overline{g}, \xi, \overline{S}}, \overline{Z}^{T, \overline{g}, \xi, \overline{S}}, \overline{K}^{T, \overline{g}, \xi, \overline{S}})$
 of RBSDE (2.2) with data $(\xi, \overline{g}, \overline{S}) $ on the interval
$[0, T]$.
  Obviously,  $\varepsilon^{r, S}_{g,\tau}[\xi] = Y^{\tau, g, \xi,
  S}_0$ and
  $\varepsilon^{r, \overline{S}}_{\overline{g},T}[\xi] = \overline{Y}^{T,\overline{g}, \xi, \overline{S}}_{0}.$

  For simplicity, denote $(Y^{\tau, g, \xi, S}, Z^{\tau, g, \xi, S}, K^{\tau, g, \xi,
 S})$ and $(\overline{Y}^{T, \overline{g}, \xi, \overline{S}}, \overline{Z}^{T, \overline{g}, \xi, \overline{S}}, \overline{K}^{T, \overline{g}, \xi,
 \overline{S}})$ by $(Y, Z, K)$ and $(\overline{Y}, \overline{Z},
 \overline{K})$, respectively. From Lemma 2.5, we have
$$
\begin{array}{ll}
&\overline{Y}(t) = \xi, \overline{Z}(t) = 0, \overline{K}(t) =
\overline{K}(\tau)\ \mbox{on the interval}\ (\tau, T];\\
&\overline{Y}(t) = Y(t), \overline{Z}(t) = Z(t), \overline{K}(t) =
K(t)\ \mbox{on the interval}\ [0, \tau].\\
\end {array}
$$
Therefore, \ $$\varepsilon^{r, S}_{g,\tau}[\xi] = \varepsilon^{r,
\overline{S}}_{\overline{g},T}[\xi].$$ \endpf

\bt \label{T4} Assume that two functions $g_{1}$ and $g_{2}$
satisfy assumptions (A1), (A3) and (A4), and the obstacle process
$\{S_t\}_{0\leq t \leq T}$ satisfies (A5). Moreover, assume that
there is  a constant $C$ such that
 $$\sup_{0\le t\le T} S_t \leq C \q   a.s.. $$
  If for each stopping time $\tau \leq T$,\ we
have $$\varepsilon^{r, S}_{g_{1}, \tau}[\xi] \geq \varepsilon^{r,
S}_{g_{2}, \tau}[\xi] \hbox{ \it for any } \xi \in L^{2}(\Omega,
{\cal{F}}_{\tau}, P) \hbox{ \it such that } \xi  \geq S_{\tau}\
a.s.,$$ then we have
$$g_{1}(t,y,z)\geq g_{2}(t,y,z) \q a.s. \hbox{ \it for any } (t, y,
z)\in [0,T]\times [C, +\infty)\times \mathbf{R}^d.$$\et

{\bf Proof.} For each $\delta > 0$ and  $(y, z)\in
(C,+\infty)\times {\mathbf{R}}^{d}$, define the following stopping
time:
$$\tau_{\delta}=\tau_{\delta}(y,z)= \mbox{inf}\{ t\geq 0: g_{1}(t,y,z) \leq g_{2}(t,y,z)-\delta\}\wedge T. $$
If the result does not hold, then there exists $\delta > 0 $ and
$(y, z)\in (C,+\infty)\times {\mathbf{R}}^{d}$ such that
$$P(\{\tau_{\delta}(y,z) < T\}) > 0.$$

For such a triple $(\delta, y, z),$ consider the following SDEs
defined on the interval $[\tau_{\delta}, T]$:
$$
  \left\{
  \begin{array}{ll}
   -dY^{1}(t) = g_{1}(t, Y^{1}(t), z)dt  - z dB_{t},\\
   Y^{1}(\tau_{\delta})= y\\
  \end{array}
  \right.
   $$
   and
   $$
  \left\{
  \begin{array}{ll}
   -dY^{2}(t) = g_{2}(t, Y^{2}(t), z)dt  - z dB_{t},\\
   Y^{2}(\tau_{\delta})= y.
  \end{array}
  \right.
   $$
    For $i = 1, 2,$\ the above equations admit a unique solution
$Y^{i} \in {\cal{S}}^2(\tau_{\delta}, T; {\mathbf{R}})$.

 Now we define  the following stopping times:
$$\tau^1_{\delta}= \mbox{inf}\{ t\geq \tau_{\delta}: Y^{1}_t\leq S_t\}\wedge T, $$
$$\tau^2_{\delta}= \mbox{inf}\{ t\geq \tau_{\delta}: Y^{2}_t\leq S_t\}\wedge T, $$
$$\tau'_{\delta}= \mbox{inf}\{ t\geq \tau_{\delta}: g_{1}(t,Y^{1}(t),z) \geq g_{2}(t,Y^{2}(t),z)-\frac{\delta}{2}\}\wedge T. $$
Note that $\tau^1_{\delta}= \tau^2_{\delta}= \tau'_{\delta}= T,$\
if $\tau_{\delta}= T$. Obviously, $\{\tau_{\delta} <
\tau^1_{\delta}\}=\{\tau_{\delta} <
\tau^2_{\delta}\}=\{\tau_{\delta} < \tau'_{\delta}\}=\{
\tau_{\delta}< T\}.  $ Define $$\tau^3_{\delta}=
\tau^1_{\delta}\wedge\tau^2_{\delta}\wedge\tau'_{\delta}. $$\
Hence $P(\{\tau_{\delta} < \tau^3_{\delta}\})
> 0.$ Moreover,  we have
$Y^1_t>S_t$ and $Y^2_t>S_t$ on the interval $[\tau_{\delta},
\tau^3_{\delta})$. Then the solution $(Y^{i}(t), z, 0)$ is  the
solution of RSBDE (2.2) with data $(Y^i(\tau^3_{\delta}), g_{i},
S)$ on the interval $[\tau_{\delta}, \tau^3_{\delta}]$ for $
i=1,2$.

    We first get the following three lemmas.
\bl $\varepsilon^{r, S}_{g_{i}, \tau_{\delta}}[y]
 = \varepsilon_{g_{i}, \tau_{\delta}}[y]= y$ for $i=1, 2.$\el

 {\bf Proof.} Consider the following BSDE defined on the interval $[0, \tau_{\delta}]$:
$$
  \left\{
  \begin{array}{ll}
   -dy^{1}(t) = g_{1}(t, y^{1}(t), z^{1}(t))dt - z^{1}(t)dB_{t},\\
   y^{1}(\tau_{\delta})= y.\\
  \end{array}
  \right.\\
   $$
   From assumption (A3), we see that
   $$y^{1}(t) =  y \hbox{ \rm and } z^{1}(t)=0 \mbox{ \rm on the interval }  [0, \tau_{\delta}]. $$
   Obviously,  the triple $(y^{1}, z^{1}, 0)$ is  the solution
   of RSBDE (2.2) with data $(y, g_{1}, S)$ on the interval
   $[0, \tau_{\delta}]$. Similarly, we have $\varepsilon^{r, S}_{g_{2}, \tau_{\delta}}[y]
 = \varepsilon_{g_{2}, \tau_{\delta}}[y]= y.$ The proof is
 complete.\endpf

\bl The strict inequality $ Y^{1}(\tau^3_{\delta}) >
Y^{2}(\tau^3_{\delta})$ holds on $\{\tau_{\delta} <
 \tau^3_{\delta}\}.$\el

{\bf Proof.}\ From the definitions of $\tau'_{\delta}$ \ and $
Y^{i}$, we have \par
 $$Y^{1}(\tau^3_{\delta}) - Y^{2}(\tau^3_{\delta})=
 \int^{\tau^3_{\delta}}_{\tau_{\delta}}[g_{2}(s,Y^{2}(s),z)-g_{1}(s,Y^{1}(s),z)]ds
  \geq \frac{\delta}{2}(\tau^3_{\delta}-\tau_{\delta})> 0,$$ $\mbox{on}\ \  \{\tau_{\delta} <
 \tau^3_{\delta}\}. $ The proof is complete.\endpf

\bl $\varepsilon^{r, S}_{g_{2},
\tau^3_{\delta}}[Y^{1}(\tau^3_{\delta})]
 = \varepsilon_{g_{2}, \tau^3_{\delta}}[Y^{1}(\tau^3_{\delta})].
 $\el

{\bf Proof.} Consider the following BSDE:
$$
  \left\{
  \begin{array}{ll}
   -d\tilde Y^{2}(t) = g_{2}(t, \tilde Y^{2}(t), \tilde Z^{2}(t))dt - \tilde Z^{2}(t)dB_{t},\q t\in [0, \tau^3_{\delta}];\\
   \tilde Y^{2}(\tau^3_{\delta})= Y^{1}(\tau^3_{\delta}).\\
  \end{array}
  \right.\\
   $$
From the definition of $\tau^3_{\delta}$ and Lemma 5.2, we get
$$Y^{1}(\tau^3_{\delta}) \geq Y^{2}(\tau^3_{\delta})
\mbox{ \rm and } P(\{Y^{1}(\tau^3_{\delta}) >
Y^{2}(\tau^3_{\delta})\})>0.
$$
On the other hand, we have
 $$\varepsilon^{r, S}_{g_{2}, \tau^3_{\delta}}[Y^{2}(\tau^3_{\delta})]
 = \varepsilon_{g_{2}, \tau^3_{\delta}}[Y^{2}(\tau^3_{\delta})] \hbox{ \rm and } \varepsilon^{r,
S}_{g_{2},\tau^3_{\delta}}[Y^{2}(\tau^3_{\delta})|{\cal{F}}_t]=\varepsilon_{g_{2},
\tau^3_{\delta}}[Y^{2}(\tau^3_{\delta})|{\cal{F}}_t]\ \ \mbox{on}\
\ [0, \tau^3_{\delta}]. $$
 From Lemma 2.1, we get
$\tilde Y^{2}(t)\geq\varepsilon_{g_{2},
\tau^3_{\delta}}[Y^{2}(\tau^3_{\delta})|{\cal{F}}_t]\geq S(t)\
\mbox{on}\ \ [0, \tau^3_{\delta}].$ Therefore $(\tilde Y^{2},
\tilde Z^{2}, 0)$ is the solution of RBSDE (2.2) with data
$(Y^1(\tau^3_{\delta}), g_{2}, S)$ on the interval $[0,
\tau^3_{\delta}]$. The proof is complete.\endpf

 Let us return to the proof of Theorem 5.1.

 Thanks to Lemma 5.1, we have
 $$y =\varepsilon^{r, S}_{g_{1}, \tau_{\delta}}[y] =
 \varepsilon^{r, S}_{g_{1}, \tau_{\delta}}[\varepsilon^{r, S}_{g_{1},
 \tau^3_{\delta}}[Y^{1}(\tau^3_{\delta})|{\cal{F}}_{\tau_{\delta}}]]
 = \varepsilon^{r, S}_{g_{1}, \tau^3_{\delta}}[Y^{1}(\tau^3_{\delta})],$$
and
$$y =\varepsilon^{r, S}_{g_{2}, \tau_{\delta}}[y] =
 \varepsilon^{r, S}_{g_{2}, \tau_{\delta}}[\varepsilon^{r, S}_{g_{2},
 \tau^3_{\delta}}[Y^{2}(\tau^3_{\delta})|{\cal{F}}_{\tau_{\delta}}]]
 = \varepsilon^{r, S}_{g_{2}, \tau^3_{\delta}}[Y^{2}(\tau^3_{\delta})].$$

On the other hand, from the definition of $\tau_\delta^3$ and
Lemma 5.3, it follows that
 $$\varepsilon^{r, S}_{g_{2}, \tau^3_{\delta}}[Y^{2}(\tau^3_{\delta})]
 = \varepsilon_{g_{2}, \tau^3_{\delta}}[Y^{2}(\tau^3_{\delta})]
 $$
 and
$$\varepsilon^{r, S}_{g_{2}, \tau^3_{\delta}}[Y^{1}(\tau^3_{\delta})]
 = \varepsilon_{g_{2}, \tau^3_{\delta}}[Y^{1}(\tau^3_{\delta})],
 $$
 respectively. Furthermore, from Lemma 2.2 we get
 $$ \varepsilon^{r, S}_{g_{2}, \tau^3_{\delta}}[Y^{i}(\tau^3_{\delta})]= \varepsilon_{g_{2}, \tau^3_{\delta}}[Y^{i}(\tau^3_{\delta})]=
 \varepsilon_{\overline{g}_2, T}[Y^{i}(\tau^3_{\delta})], \ i = 1, 2. $$
Here, $\overline{g}_2(t, y, z) := g_2(t, y, z)1_{[0,
\tau^3_{\delta}]}(t)$ for $\ae t \in [0, T]$ and any $(y, z)\in
{\mathbf{R}}\times {\mathbf{R}}^d.$ From the definition of
$\tau^3_{\delta}$ and Lemma 5.2, it follows that
$$Y^{1}(\tau^3_{\delta}) \geq Y^{2}(\tau^3_{\delta})
\mbox{ \rm and } P(\{Y^{1}(\tau^3_{\delta}) >
Y^{2}(\tau^3_{\delta})\})>0.
$$ Therefore, in view of Lemma~\ref{L1}, we have
 \be \varepsilon_{\overline{g}_2,T}[Y^2(\tau_\d^3)]<\varepsilon_{\overline{g}_2,T}[Y^1(\tau_\d^3)]. \ee
 Concluding the above, we get
$$
\begin{array}{ll}
y &= \varepsilon^{r, S}_{g_{2},
\tau^3_{\delta}}[Y^{2}(\tau^3_{\delta})]
 = \varepsilon_{g_{2}, \tau^3_{\delta}}[Y^{2}(\tau^3_{\delta})]=\varepsilon_{\overline{g}_2, T}[Y^{2}(\tau^3_{\delta})]
 < \varepsilon_{\overline{g}_2, T}[Y^{1}(\tau^3_{\delta})]=\varepsilon_{g_{2},
 \tau^3_{\delta}}[Y^{1}(\tau^3_{\delta})]\\
 &= \varepsilon^{r, S}_{g_{2}, \tau^3_{\delta}}[Y^{1}(\tau^3_{\delta})]
\leq \varepsilon^{r, S}_{g_{1},
\tau^3_{\delta}}[Y^{1}(\tau^3_{\delta})]=
y.\\
\end{array}
$$ This is a contradiction. The proof is complete. \endpf

\br Consider the example given in Remark 4.2. Furthermore, assume
that $\mu_2>0$. Immediately, we have the following three facts:
(i) $g_1(\cdot, y, \cdot)=g_2(\cdot, y, \cdot)$ when $y\geq c_2$;
(ii) $g_1(\cdot, y, z)<g_2(\cdot, y, z)$ when $c_1\leq y<c_2$ and
$z\neq 0$; and (iii) $g_1(\cdot, y, z)>g_2(\cdot, y, z)$ when
$y\leq c_1- |z|$ and $z\neq 0$.

On the other hand, since
$$\varepsilon^{r}_{g_{1}, T}[\xi] = \varepsilon^{r}_{g_{2},
T}[\xi] \hbox{ \it and}\
 \varepsilon^{r}_{g_{1}, T}[\xi|{\cal{F}}_{t}] = \varepsilon^{r}_{g_{2}, T}[\xi|{\cal{F}}_{t}]\ \ a.s. \hbox{ \it for } t\in (0,
    T]$$
    for any $\xi \in L^{2}(\Omega,
{\cal{F}}_{T}, P)$ satisfying $\xi \geq S_{T}\ a.s.$, we deduce
the above fact (i) from Theorem 5.1. The other two facts (ii) and
(iii) demonstrate that the conclusion of Theorem 5.1 is the best
possible in the underlying example. \er

In Theorem 5.1, the bound assumption on the obstacle process
appears to be very restrictive. In what follows, we show that if
the generator of RBSDE~(\ref{RBSDE}) does not depend on the first
unknown variable $y$, we can get the following global converse
comparison result without the bound assumption.

\bt \label{T5} Suppose that two fields $g_{1}$ and $g_{2}$ satisfy
assumptions (A1), (A3) and (A4), and the obstacle process
$\{S_t\}_{0\leq t \leq T}$ satisfies (A5). Furthermore, assume
that $g_{1}$ and $g_{2}$ do not depend on $y$.  If for each
stopping time $\tau \leq T$, \be \varepsilon^{r, S}_{g_{1},
\tau}[\xi] \geq \varepsilon^{r, S}_{g_{2}, \tau}[\xi] \hbox{ \it
for any } \xi \in L^{2}(\Omega, {\cal{F}}_{\tau}, P) \hbox{ \it
such that } \xi \geq S_{\tau}\ a.s.,\label{1*}\ee then we have \be
g_{1}(t,z)\geq g_{2}(t,z) \q \as \hbox{ \it for } (t, z)\in [0,
T]\times{\mathbf{R}^d}.\label{2*}\ee \et

{\bf Proof. } {\bf Step 1.}  If $\sup_{0\le t\le T}S_t$ is bounded
from above, then the desired assertion is immediate.

{\bf Step 2. } For a large integer $n$, define the stopping time
$$\tau_n := \mbox{inf}\{t\geq 0: S_t \geq n \}\wedge T.$$
Then  $0 \leq \tau_n \leq T \, \as$.   Since $S_0$\ is a
deterministic finite number and $S$ is continuous, we have
$\tau_n>0 \ a.s.$ for any $n> S_0+1.$

For every $n > S_0+1,$  define $\overline{g}_{i}(t, z):=
g_{i}(t,z)1_{[0, \tau_{n}]}(t)$ and $\overline{S}_{t}=S_{t\wedge
\tau_{n}}$ for $ (t, z)\in [0, T]\times{\mathbf{R}^d}$ with $i=1,
2.$ Then for each stopping time $\tau \leq T$ and any $\xi \in
L^{2}(\Omega, {\cal{F}}_{\tau}, P)$ such that $\xi  \geq
\overline{S}_{\tau}$, if we have \be \varepsilon^{r,
\overline{S}}_{\overline{g}_{1}, \tau}[\xi] \geq \varepsilon^{r,
\overline{S}}_{\overline{g}_{2}, \tau}[\xi],\ \label{*}  \ee then
noting that  $\overline{S}_t \leq n$ on $[0, T]$ (in view of  the
definition of $\tau_{n}$),
 we have from Step 1 that
  $$\overline{g}_{1}(t,z)\geq \overline{g}_{2}(t,z)\q \as \hbox{ \rm for } (t, z)\in [0,
T]\times{\mathbf{R}^d}.$$ That is,
$$g_{1}(t,z)\geq g_{2}(t,z)\q \as \hbox{ \rm for } (t, z)\in [0,
\tau_n]\times{\mathbf{R}^d}.$$ Obviously, $\tau_{n}\uparrow T$ as
$n\rightarrow\infty.$ Passing to limit, from assumption (A4) we
get $$g_{1}(t,z)\geq g_{2}(t,z)\q \as \hbox{ \rm for } (t, z)\in
[0, T]\times{\mathbf{R}^d}.$$ The proof is then complete.
Therefore we only need to prove inequality~(\ref{*}).

 Define
${\tilde{g}}_{i}(t, z)= \overline{g}_{i}(t,z)1_{[0, \tau]}(t)$ and
${\tilde{S}}_{t}=\overline{S}_{t\wedge \tau}$ for $(t, z)\in [0,
T]\times{\mathbf{R}^d}$ with $i=1, 2.$ It follows that
${\tilde{g}}_{i}(t, z)= g_{i}(t,z)1_{[0, \tau\wedge\tau_{n}]}(t)$
and ${\tilde{S}}_{t}=S_{t\wedge \tau\wedge\tau_{n}}$ for $(t,
z)\in [0, T]\times{\mathbf{R}^d},\ i=1, 2.$ From Proposition 5.1,
we have
$$\varepsilon^{r, \overline{S}}_{\overline{g}_{i}, \tau}[\xi]=\varepsilon^{r, \tilde{S}}_{\tilde{g}_{i}, T}[\xi], \ i=1, 2.$$

On the other hand, from the definitions of ${\tilde{g}}_{1}(t,
z)$\ and ${\tilde{S}}$,\ we have
$$\varepsilon^{r, \tilde{S}}_{\tilde{g}_{1}, T}[\xi]=
\varepsilon^{r, \tilde{S}}_{\tilde{g}_{1},
\tau\wedge\tau_{n}}[\varepsilon^{r, \tilde{S}}_{\tilde{g}_{1},
T}[\xi|{\cal{F}}_{\tau\wedge\tau_{n}}]]= \varepsilon^{r,
S}_{g_{1}, \tau\wedge\tau_{n}}[\varepsilon^{r, \tilde{S}}_{\tilde
{g}_{1}, T}[\xi|{\cal{F}}_{\tau\wedge\tau_{n}}]].$$ Therefore
$$\varepsilon^{r, \overline{S}}_{\overline{g}_{1}, \tau}[\xi]=\varepsilon^{r,
S}_{g_{1}, \tau\wedge\tau_{n}}[\varepsilon^{r, \tilde{S}}_{\tilde
{g}_{1}, T}[\xi|{\cal{F}}_{\tau\wedge\tau_{n}}]].$$ Similarly,
$$\varepsilon^{r, \overline{S}}_{\overline{g}_{2}, \tau}[\xi]=\varepsilon^{r,
S}_{g_{2}, \tau\wedge\tau_{n}}[\varepsilon^{r, \tilde {S}}_{\tilde
{g}_{2}, T}[\xi|{\cal{F}}_{\tau\wedge\tau_{n}}]].$$
 Also, thanks to the
definitions of ${\tilde {g}}_{1}(t, z)$,\ ${\tilde{g}}_{2}(t, z)$\
and ${\tilde {S}}$, we get
$$\varepsilon^{r,
\tilde {S}}_{\tilde{g}_{1},
T}[\xi|{\cal{F}}_{\tau\wedge\tau_{n}}]=\varepsilon^{r,
\tilde{S}}_{\tilde{g}_{2},
T}[\xi|{\cal{F}}_{\tau\wedge\tau_{n}}].$$ For simplicity,  we set
$\eta :=\varepsilon^{r, \tilde{S}}_{\tilde{g}_{1},
T}[\xi|{\cal{F}}_{\tau\wedge\tau_{n}}].$ Obviously $\eta \in
L^{2}(\Omega, {\cal{F}}_{\tau\wedge\tau_{n}}, P)$ and $ \eta  \geq
S_{\tau\wedge\tau_{n}}.$ Then from the assumption, it follows that
$$\varepsilon^{r, \overline{S}}_{\overline{g}_{1}, \tau}[\xi]=\varepsilon^{r,
S}_{g_{1}, \tau\wedge\tau_{n}}[\eta]\geq \varepsilon^{r,
S}_{g_{2}, \tau\wedge\tau_{n}}[\eta]=\varepsilon^{r,
\overline{S}}_{\overline{g}_{2}, \tau}[\xi].$$ Now we end up with
the proof. \endpf

\br Obviously, (\ref{1*}) and (\ref{2*}) in Theorem 5.2
 are also equivalent.\er

If assumption (A3) is replaced with assumption (A2) in
Theorem~\ref{T4}, then we have

 \bt \label{T6} Assume that two random fields $g_{1}$ and $g_{2}$
satisfy assumptions (A1), (A2) and (A4), and the obstacle process
$\{S_t\}_{0\leq t \leq T}$ satisfies (A5).
  If for any two stopping times
$\tau$ and $\sigma$ such that  $\tau \leq\sigma\leq T$,
  \be\varepsilon^{r, S}_{g_{1}, \sigma}[\xi|{\cal{F}}_{\tau}] \geq
\varepsilon^{r, S}_{g_{2}, \sigma}[\xi|{\cal{F}}_{\tau}]\ a.s.
\hbox{ \it for } \xi \in L^{2}(\Omega, {\cal{F}}_{\sigma}, P)
\hbox{ \it such that } \xi \geq S_{\sigma}\ a.s.,\ee then for  any
continuous process $Y\in {\cal{S}}^2(0, T; {\mathbf{R}})$ such
that $Y(t)\geq S_t\ a.s.$  with $t\in [0,T]$, we have \be
g_{1}(t,Y(t),z)\geq g_{2}(t,Y(t),z)\q a.s. \hbox{ \it for } (t,
 z)\in [0,T]\times\mathbf{R}^d.\label{last}\ee
 In particular, \be g_{1}(t,S(t),z)\geq g_{2}(t,S(t),z)\q a.s.
\hbox{ \it for } (t,
 z)\in [0,T]\times\mathbf{R}^d.\ee\et

{\bf Proof.} In view of the continuity of $g_1(t,y,z)$ and
$g_2(t,y,z)$ in $y$, it is sufficient to prove~(\ref{last}) for
any continuous process $Y\in {\cal{S}}^2(0, T; {\mathbf{R}})$ such
that $Y(t)\geq S_t + \epsilon\ a.s.$  with $t\in [0,T]$ for some
constant $\epsilon>0$. We shall prove it by contradiction.

Otherwise, there would exist $\delta
> 0 $ and $z\in {\mathbf{R}}^{d}$ such that
$$P(\{\tau_{\delta}(z) < T\}) > 0.$$
Here for $\delta > 0 \ \mbox{and}\ z\in {\mathbf{R}}^{d}$, we have
defined the following stopping time:
$$\tau_{\delta}=\tau_{\delta}(z):= \mbox{inf}\{ t\geq 0: g_{1}(t,Y(t),z) \leq g_{2}(t,Y(t),z)-\delta\}\wedge T. $$

For such a pair $(\delta, z),$ analogous to the proof of
Theorem~\ref{T4}, consider the following SDEs defined on the
interval $[\tau_{\delta}, T]$:
$$
  \left\{
  \begin{array}{ll}
   -dY^{1}(t) = g_{1}(t, Y^{1}(t), z)dt  - z dB_{t},\\
   Y^{1}(\tau_{\delta})= Y(\tau_\d)\\
  \end{array}
  \right.\\
   $$
   and
   $$
  \left\{
  \begin{array}{ll}
   -dY^{2}(t) = g_{2}(t, Y^{2}(t), z)dt - z dB_{t},\\
   Y^{2}(\tau_{\delta})= Y(\tau_\d).\\
  \end{array}
  \right.\\
   $$
  The above SDEs admit unique solutions
   $Y^{i} \in {\cal{S}}^2(\tau_{\delta}, T; {\mathbf{R}})$ with $i = 1,
   2$.

Define  the following stopping times:
$$\tau^1_{\delta}= \mbox{inf}\{ t\geq \tau_{\delta}: Y^{1}_t\leq S_t\}\wedge T, $$
$$\tau^2_{\delta}= \mbox{inf}\{ t\geq \tau_{\delta}: Y^{2}_t\leq S_t\}\wedge T,
$$ and
$$\tau'_{\delta}= \mbox{inf}\{ t\geq \tau_{\delta}: g_{1}(t,Y^{1}(t),z) \geq g_{2}(t,Y^{2}(t),z)-\frac{\delta}{2}\}\wedge T. $$
Note that $\tau^1_{\delta}= \tau^2_{\delta}= \tau'_{\delta}= T,$\
if $\tau_{\delta}= T$. Obviously, $\{\tau_{\delta} <
\tau^1_{\delta}\}=\{\tau_{\delta} <
\tau^2_{\delta}\}=\{\tau_{\delta} < \tau'_{\delta}\}=\{
\tau_{\delta}< T\}.  $ We define $$\tau^3_{\delta}=
\tau^1_{\delta}\wedge\tau^2_{\delta}\wedge\tau'_{\delta}. $$ Hence
$P(\{\tau_{\delta} < \tau^3_{\delta}\})
> 0.$ Moreover, we have
$Y^1_t>S_t$ and $Y^2_t>S_t$\ on the interval $[\tau_{\delta},
\tau^3_{\delta})$. Therefore, the triple $(Y^{i}, z, 0)$ is the
solution of RSBDE (2.2) with data $(Y^i(\tau^3_{\delta}), g_{i},
S)$ on the interval $[\tau_{\delta}, \tau^3_{\delta}]$ for
$i=1,2$. Consequently,
 $$\varepsilon^{r,
S}_{g_{1},\tau^3_{\delta}}[Y^{1}(\tau^3_{\delta})|{\cal{F}}_{\tau_{\delta}}]
=\varepsilon_{g_{1},\tau^3_{\delta}}[Y^{1}(\tau^3_{\delta})|{\cal{F}}_{\tau_{\delta}}]=Y(\tau_\d)$$
and
$$\varepsilon^{r,
S}_{g_{2},\tau^3_{\delta}}[Y^{2}(\tau^3_{\delta})|{\cal{F}}_{\tau_{\delta}}]
=\varepsilon_{g_{2},\tau^3_{\delta}}[Y^{2}(\tau^3_{\delta})|{\cal{F}}_{\tau_{\delta}}]=Y(\tau_\d).$$

Identical to the proof of lemmas 5.2 and 5.3, we get

 \bl We have
\label{L4}\be Y^{1}(\tau^3_{\delta})
> Y^{2}(\tau^3_{\delta})\  \hbox{ \it on } \{\tau_{\delta} <
 \tau^3_{\delta}\}\label{ine}\ee
 and
\be \label{equ}\varepsilon^{r,
S}_{g_{2},\tau^3_{\delta}}[Y^{1}(\tau^3_{\delta})|{\cal{F}}_{\tau_{\delta}}]
=\varepsilon_{g_{2},\tau^3_{\delta}}[Y^{1}(\tau^3_{\delta})|{\cal{F}}_{\tau_{\delta}}].\ee
\el

From the definition of $\tau^3_{\delta}$ and~(\ref{ine}), we have
$$Y^{1}(\tau^3_{\delta}) \geq Y^{2}(\tau^3_{\delta})\ a.s.
\mbox{ \rm and } P(\{Y^{1}(\tau^3_{\delta}\}) >
Y^{2}(\tau^3_{\delta}))>0.
$$
Then it follows from Lemmas~\ref{L1} and~(\ref{equ}) that
$$\varepsilon^{r,
S}_{g_{2},\tau^3_{\delta}}[Y^{1}(\tau^3_{\delta})|{\cal{F}}_{\tau_{\delta}}]
=\varepsilon_{g_{2},\tau^3_{\delta}}[Y^{1}(\tau^3_{\delta})|{\cal{F}}_{\tau_{\delta}}]\geq
\varepsilon_{g_{2},\tau^3_{\delta}}[Y^{2}(\tau^3_{\delta})|{\cal{F}}_{\tau_{\delta}}]=\varepsilon^{r,
S}_{g_{2},\tau^3_{\delta}}[Y^{2}(\tau^3_{\delta})|{\cal{F}}_{\tau_{\delta}}]\
a.s.
$$
and
$$P(\{\varepsilon^{r,
S}_{g_{2},\tau^3_{\delta}}[Y^{1}(\tau^3_{\delta})|{\cal{F}}_{\tau_{\delta}}]
>\varepsilon^{r,
S}_{g_{2},\tau^3_{\delta}}[Y^{2}(\tau^3_{\delta})|{\cal{F}}_{\tau_{\delta}}]\})=P(\{\varepsilon_{g_{2},\tau^3_{\delta}}[Y^{1}(\tau^3_{\delta})|{\cal{F}}_{\tau_{\delta}}]
>\varepsilon_{g_{2},\tau^3_{\delta}}[Y^{2}(\tau^3_{\delta})|{\cal{F}}_{\tau_{\delta}}]\})>0.
$$
The last relation implies that
$$P(\{\varepsilon^{r,
S}_{g_{2},\tau^3_{\delta}}[Y^{1}(\tau^3_{\delta})|{\cal{F}}_{\tau_{\delta}}]
>Y(\tau_\d)\})>0$$
which contradicts the assumption that
$$
\varepsilon^{r,
S}_{g_{2},\tau^3_{\delta}}[Y^{1}(\tau^3_{\delta})|{\cal{F}}_{\tau_{\delta}}]
\leq\varepsilon^{r,
S}_{g_{1},\tau^3_{\delta}}[Y^{1}(\tau^3_{\delta})|{\cal{F}}_{\tau_{\delta}}]=Y(\tau_\d)\q
a.s..
$$  The proof is complete. \endpf

The following gives an immediate consequence of Theorem~\ref{T6}.

  \bc  Suppose that two generators
$g_{1}, g_{2}$ satisfy assumptions (A1), (A2) and (A4), and the
obstacle process $\{S_t\}_{0\leq t \leq T}$ satisfies (A5).
Furthermore, assume that $g_{1}$ and $g_{2}$\ do not depend on the
first unknown variable $y$. If for each pair of stopping times
$\tau$ and $\sigma$ such that $\tau \leq\sigma\leq T$,\ we have
$$\varepsilon^{r, S}_{g_{1}, \sigma}[\xi|{\cal{F}}_{\tau}] \geq
\varepsilon^{r, S}_{g_{2}, \sigma}[\xi|{\cal{F}}_{\tau}]\ a.s.
\hbox{ \it for any } \xi \in L^{2}(\Omega, {\cal{F}}_{\sigma}, P)
\hbox{ \it such that } \xi \geq S_{\sigma}\ a.s.,$$ then we have
$$g_{1}(t, z)\geq g_{2}(t, z) \q  a.s. \hbox{ \it for any } (t, z)\in [0,
T]\times{\mathbf{R}^d}.$$\ec

\noindent {\bf Acknowledgements.}
  Both authors would like to thank Professor Shige Peng  for his helpful
  comments, and  Guangyan Jia and Zhiyong Yu for
  detecting an error in the original proof of Theorem 3.1. They
  are also very grateful to the Associate Editor and the anonymous
  referee for their helpful suggestions and criticisms.

\end{document}